\documentclass[12pt]{amsart} 
\newfont{\sheaf}{eusm10 scaled\magstep1}

\newcommand{\ra}{\ensuremath{\rightarrow}}

\def\eea{\end{eqnarray*}}
\def\bea{\begin{eqnarray*}}

\newcommand{\Proof}{{\it Proof. }} 
\newcommand{\QED}{{\hfill $Q.E.D.$}} 
 
\newtheorem{teo}{Theorem}[section] 
\newtheorem{df}[teo]{Definition} 
\newtheorem{lem}[teo]{Lemma} 
\newtheorem{cor}[teo]{Corollary} 
\newtheorem{ex}[teo]{Example} 
\newtheorem{oss}[teo]{Remark}

\newcommand{\C}{\ensuremath{\mathbb{C}}}
\newcommand{\R}{\ensuremath{\mathbb{R}}} 
\newcommand{\Z}{\ensuremath{\mathbb{Z}}} 
\newcommand{\Q}{\ensuremath{\mathbb{Q}}}
\newcommand{\F}{\ensuremath{\mathbb{F}}}

\newcommand{\HH}{\ensuremath{\mathbb{H}}} 
\newcommand{\PP}{\ensuremath{\mathbb{P}}}

\newcommand{\FF}{\ensuremath{\mathcal{F}}}
\newcommand{\EE}{\ensuremath{\mathcal{E}}}

\begin{document}

\title{Fibred K\"ahler and quasiprojective groups.}

\author{Fabrizio Catanese\\
(Universit\"at Bayreuth)\\
  This article is dedicated to Adriano Barlotti on 
the occasion of his 80-th birthday.
}

\footnote{
The research of the  author was performed in the realm  of the 
 SCHWERPUNKT "Globale Methode in der komplexen Geometrie",
and of the EAGER EEC Project. 
}
 \\

\date{February 16, 2003}
\maketitle

\begin{abstract} 
 We formulate a new  theorem giving  several necessary 
and sufficient conditions in order that a surjection 
of the fundamental group $\pi_1(X)$ of a compact
K\"ahler manifold onto the fundamental group $\Pi_g$ of a compact
 Riemann surface of genus $g \geq 2$ be induced by a holomorphic
map. For instance, it suffices that the kernel be finitely
generated.

We derive as a corollary a restriction for a group 
$G$, fitting into an exact sequence $ 1 \ra H \ra G \ra \Pi_g \ra 1$,
where $H$ is finitely generated,
 to be the fundamental group of a compact K\"ahler manifold. 

Thanks to  the extension by Bauer and Arapura of the 
Castelnuovo de Franchis theorem  to the quasi-projective case 
(more generally, to Zariski open sets of compact
K\"ahler manifolds) we first extend the previous result to the
non compact case.
We are finally able to give a topological characterization of
 quasi-projective surfaces 
which are 
fibred over a (quasi-projective) curve by a  proper holomorphic map 
of maximal rank, and we extend the previous restriction to the 
monodromy of any fibration onto a curve.

  \end{abstract}

\vfill
\pagebreak

\section{Introduction }

The history of fibrations of algebraic  (or K\"ahler) manifolds 
$f : X \ra C$ over curves $C$
of genus at least $2$, called classicaly irrational pencils, 
has a long history.

Around 1905 almost simultaneously De Franchis and Castelnuovo-Enriques 
(\cite{df05}, \cite{cast05-2}) found
that the existence of such a fibration is equivalent to the existence
of at least $2$ linearly independent holomorphic $1$-forms whose wedge 
product yields a holomorphic $2$-form which is identically zero.

Hodge theory was yet to be developed and only much later (\cite{cat1}) it
was shown that combining the Hodge decomposition with the theorem of 
Castelnuovo De Franchis one obtains a topological characterization
of such fibrations via any subspace in De Rham cohomology obtained as
the pull-back of a maximal isotropic subspace in the cohomology of $C$.

Other topological characterizations in terms of the induced surjection
of fundamental groups $ f_* : \pi_1(X) \ra \pi_1(C)$, or other statements 
in this direction, had been earlier 
obtained by several authors (cf. Jost-Yau and Siu, \cite{siu1} ,
\cite{j-y1}, \cite{j-y2}, \cite{siu2}, 
who used the theory of harmonic maps,  Beauville, \cite{bea}, used instead
 the generic 
vanishing theorems of Green and Lazarsfeld, \cite{g-l1}).

In \cite{cat2} I tried to show how the isotropic subspace theorem,
which predicts the genus of the image curve $C$ (it equals to the dimension
of the corresponding maximal isotropic subspace), unlike the other statements,
could be used to obtain also simple proofs of statements concerning surjections
 of fundamental groups, as the one given by Gromov (\cite{gro1}).

Hodge theory also works for quasi-projective manifolds, and it turns out
that in the non compact case it works much better than  the other methods
(\cite{ingrid}, \cite{ara}). These results were then used in \cite{cat3} 
to give
topological characterizations 
 of varieties isogenous to a product and of isotrivial
 fibrations of surfaces. 

Kotschick (\cite{kot}) instead used very similar methods to give a topological 
characterization of Kodaira fibrations, which was independently also obtained by Hillman
(\cite{hil}).

Our first motivation was to extend this result to any fibration, a goal that
we have achieved in the following 

{\bf Theorem 6.4.  } {\em Assume that $U$ is a non complete Zariski open set 
of an algebraic surface
 and that the following properties hold:

\begin{itemize}
\item
(P1) we have an exact sequence  
$$ 1 \ra \Pi_r \ra \pi_1(U) \ra \F_g \ra 1$$
where $g \geq 2$. 
\item
(P2) The topological 
Euler-Poincare' characteristic of $U$,
$e(U)$, equals $ 2 (g-1) (r-1)$.
\item
(P3) For each end $\EE$ of $U$, the corresponding  fundamental group 
$\pi_1^{\EE}$ surjects onto a cyclic subgroup of $\F_g$, and each
simple geometric generator $\gamma_i$ has a non trivial image
in $\F_g$.

\end{itemize}
Then $U$ is a good open set of a fibration, more precisely,
there exists a proper holomorphic submersion $f: U \ra C$ inducing 
the previous exact sequence.} 

To this purpose, we started to put together the existing results, both
in the compact and in the non- compact case, with some small addition:
we refer to Theorems 4.3 and 5.4 for full statements.

We only indicate here some new results therein:

{\bf Theorem A.  }

{\em Let $X$ be a compact K\"ahler manifold, and assume that its fundamental
group admits a non trivial homomorphism $\psi$ 
to the fundamental group $\Pi_g$ of a compact
 Riemann surface of genus $g \geq 2$, with kernel $H$.

Then the following conditions are equivalent 
\begin{itemize}

\item
4) $\psi$ is induced by an irrational pencil of genus $g$  
without multiple fibres
\item
5) $\psi$ is surjective and its kernel $H$ is finitely generated

\end{itemize} }

{\bf Theorem A'.  }

{\em 
Let $X$ be a compact K\"ahler manifold and $Y = X - D$ be a Zariski open 
set. Assume that the fundamental group of $Y$  admits a 
homomorphism $\psi : \pi_1(Y) \ra \F_g$ to a free group of rank $g$, 
with kernel $H$.

Then the following are equivalent 

\begin{itemize}
\item
1') $\psi$ is induced by a pencil $ f : Y \ra C$ of 
type $g$ without multiple fibres and by a surjection $\pi_1(C) \ra \F_g$.

\item
2') $\psi$ is surjective and the kernel $H$ of $\psi$ 
is finitely generated

\end{itemize}
}

The new ingredient here is a remarkable property of  free groups
 and of the fundamental groups of compact curves of genus $g \geq 2$. 

This property is that every  non trivial {\bf N}ormal subgroup of
{\bf I}nfinite index
is {\bf N}ot {\bf F}initely generated.

We abbreviate this property by the acronym {\bf NINF}, and we devote
section $3$ to establishing this result for the above mentioned groups.

This property plays an important role, for instance it shows that,
contrary to what is stated by some author, the kernel of the homomorphism 
 between fundamental groups $f_* : \pi_1(X) \ra \pi_1(C)$ needs not be
 finitely
generated, it is finitely generated
 if and only if there are no multiple fibres.

In the case where there are multiple fibres, one can take a ramified
base change which eliminates the multiple fibres, and indeed one could
extend Theorems A, A' to include the case where $H$ is not finitely
generated, but we believed that the results stated in this article are 
already sufficiently complicated, so we omitted  to treat this extension.

We prefered instead to concentrate on some important consequence,
concerning the monodromy  of fibrations over curves.
For instance, putting together  Theorem $A$ with the old 
isotropic subspace theorem we obtain

{\bf Corollary 7.3}

{\em If a finitely presented group $\Gamma$ admits a surjection 
$\Gamma \rightarrow 
\Pi_{g}$ with finitely generated kernel $H$, then $\Gamma$ cannot be 
the fundamental group of a compact K\"ahler manifold $X$ if there is
a non zero element $ u \in  H^1(H, \Z)^{\Pi _g }$ such that cup product 
with $u$ yields the identically 
  zero map $$  H^1(\Pi _g , \Z) \ra H^1(\Pi _g , H^1(H, \Z)).$$ }

We then spell out in detail the meaning of degeneracy of the
above cup product: it means that there exists a bad monodromy submodule.

By extending everything in the non compact case one obtains then
a restriction for the monodromy of fibrations over curves which is in the
 same spirit as Deligne's semisimplicity Theorem (4.2.6 of \cite{del71}).

\bigskip

\section{Notation}
\begin{itemize}
\item
 $\Pi_g$ denotes the fundamental group of a compact
 Riemann surface $C_g$ of genus $g \geq 2$,
$\Pi_g : = < a_1, \dots a_g,b_1, \dots b_g | 
[a_1, b_1]  \dots [a_g, b_g] =1 > .$

\item
$ \F_g$ denotes a free group of rank $g \geq 2$
\item
$X$ will be a compact K\"ahler manifold

\end{itemize}

\bigskip
\section{Non finitely generated subgroups }

This section is devoted to a remarkable property enjoyed, for 
$g \geq 2$, by
the free groups $ \F_g$ and by the  fundamental groups $\Pi_g$ 
of a compact
 Riemann surface $C_g$ of genus $g \geq 2$.

\begin{df}
A group $G$ is said to satisfy property $NINF$ (to be more precise,
 but less concise, we should call it $NIINFG$) if every 
{\bf normal } non trivial subgroup $K$ of {\bf infinite index }
is {\bf non finitely} generated. 

We shall also say that $G$ is a $NINF$.
\end{df}

The application of the  above notion that we shall need is the following

\begin{lem}
Let $ 1 \ra A \ra \Pi \ra B \ra 1$ be an exact sequence 
of group homomorphisms such that

\begin{itemize}
\item
1) $A$ is finitely generated 
\item
2) $B$ is infinite
\item
$\phi :  \Pi \ra B $ factors as $ \rho \circ \psi$, where 
$\psi : \Pi \ra G$ is surjective
\item
3) $G$ is a $NINF$

\end{itemize}

Then $\rho : G \ra B$ is an isomorphism.
\end{lem}

\Proof
Let $ j : A \ra \Pi$ be the inclusion and define $K' := ker 
( \psi \circ  j )$, $K : = ker 
( \psi)$. Then $ K' = K$ since $K \subset A = ker ( \phi)$.

Set $A' := A / K$, so $A'$ injects into $G = \Pi / K$. Moreover, 
$A'$ is normal in $G$ with quotient $B = \Pi / A$ which is infinite
by assumption. Since $A'$ is finitely generated, as a quotient of $A$,
and $G$ is a $NINF$ it follows that $A'$ is trivial.

Whence, $ A = K$ and $\rho : G \ra B$ is an isomorphism, as desired. 
\qed

\begin{lem}
A free group $\F_n$ enjoys property $NINF$.
\end{lem}

\Proof

W.l.o.g. assume we may assume $n \geq 2$. We view $\F_n$ as the 
fundamental group $\pi_1(Y)$, where $Y$ is a bouquet of $n$ circles.

Let $Z$ be the covering space corresponding to a normal subgroup $K$
of infinite index.

 $Z$ is 
indeed the  Cayley graph for the infinite group $G : = \F_n / K$
 with respect
 to the finite set of $n$ generators, $g_1, \dots g_n$, corresponding 
to the surjection $\F_n \ra G$. 

If $K$ is non trivial, then $Z$ is not simply connected,  and 
there is a non trivial minimal closed simplicial path $\xi$ based on 
the base
point $x_0=1$ : $ (x_0 = 1, x_1 = \gamma_1 , \dots x_m =
\gamma_1  \dots \gamma_m  )$, where the $\gamma_i$ 's belong to the
given set of generators $\{g_1, \dots g_n \}$.

Let $M \subset G$ be the set $ \{x_0 = 1, x_1 , \dots x_m \}$
and let $M'$ be the finite set $ M M^{-1} \subset G$.

Since $G$ is infinite , there exist infinitely many $h_{\alpha}$'s 
such that $h_{\alpha} \notin h_{\beta} M'$ for ${\alpha} \neq {\beta}$.

Whence, for ${\alpha} \neq {\beta}$ and $\forall x_i, x_j$, we have
$h_{\alpha} x_i \neq h_{\beta} x_j$. It follows that the 
cycles $ h_{\alpha} (\xi)$ are homologically independent.

A fortiori, we have shown $ Rank (H_1 (Z, \Z) = \infty$ and
$K$ is not finitely generated.

\qed

\begin{lem}
A fundamental group $\Pi_g$ enjoys property $NINF$ for $g \geq 2$.
\end{lem} 

\Proof
Let $\Gamma$ be a nontrivial normal subgroup of $\Pi: = \Pi_g$, and let 
$ f : D \ra C$ be the corresponding unramified covering of a compact 
Riemann surface of genus $g$. We have $ g \geq 2$, whence
we may view $D$ as a quotient of the upper half plane $\HH$ by the action of the 
group $\Gamma$ acting freely and properly discontinuously.

As in \cite{sieg}, thm. 4, page 35 , we consider fundamental domains, 
$\FF_{\Pi}$ resp. $\FF_{\Gamma}$, bounded by non-euclidean segments (possibly also 
 lines or half-lines).

While $\FF_{\Pi}$ has finite area, the area of $\FF_{\Gamma}$ is the area of 
$\FF_{\Pi}$ multiplied by the index of $\Gamma$, whence $\FF_{\Gamma}$
has infinite area.

Assume that $\Gamma$ be finitely generated: then 
( cf. \cite{beard} Thm. 10.1.2 , page 254) there is such a fundamental 
domain $\FF_{\Gamma}$ with finitely many sides.
Since however its area is infinite, it cannot be a non Euclidean ideal polygon,
and there are intervals in the real line $\PP^1_{\R}$ which need to be added 
for a compactification of $\FF_{\Gamma}$.

Let us recall the following standard definitions (cf.  \cite{sieg}, and especially
 \cite{se-so}, \cite{beard})

\begin{df}

\begin{itemize}
\item
1) A subgroup $\Gamma$ of $PSL(2,\R)$  acting properly discontinuosly on $\HH$ is 
called a Fuchsian group (more generally, a Fuchsian group is a conjugate in $PSL(2,\C)$
of such a subgroup).
\item 
2) $\Gamma$ is properly discontinuos if and only if it is discrete in $PSL(2,\R)$ 
\item
3) The {\bf Limit Set} $ L (\Gamma)$  is defined as 
$$L (\Gamma) : = \cup_{z \in \HH} ( \overline{ \Gamma z} \cap \PP^1_{\R}).$$ 
\item
4) Equivalently, (cf. \cite{se-so}, page 108) 

$$L (\Gamma) : = \{z  \in \PP^1_{\R}| \exists 
\gamma \in \Gamma ,\gamma  \neq 1 ,
 s.t.  \gamma (z) = z \}.$$

\item
5) A Fuchsian subgroup $\Gamma$ is said to be of the I Kind if 
$ L (\Gamma) = \PP^1_{\R} $ (else, it is said to be of the II Kind).
\end{itemize} 
\end{df}

We have (cf. e.g. Lemma 3.12.2, page 108 of \cite{se-so}) that if $\HH / \Gamma$ 
is compact, then $\Gamma$ is of the I Kind. This implies a consequence for
 our group $\Pi$.

In fact, since $\Pi$ normalizes $\Gamma$, 
then $\Pi$ carries the limit set
$L (\Gamma)$ to itself. 

Let $I$ be  an interval in the real line $\PP^1_{\R}$ which is in the boundary
 of $\FF_{\Gamma}$. Since $\Pi$ is of the I Kind, there is $ x \in I$ and 
$g \in \Pi - \{ 1\}$ such that $ g x = x$.

Since, as we observed, $ g( L (\Gamma)) = L (\Gamma)$, $g$ carries the interior of the 
complement to $L (\Gamma)$ in $\PP^1_{\R}$ to itself.  Let the interval $(a,b)$ be
the connected component of this interior containing $x$. Since moreover $ g x = x$,
$g$ carries $(a,b)$ to itself.

Assume that $ a \neq b$ : then $g^2$ has three fixed points ( $a,b,x$), thus $g^2$ is
the identity, contradicting the hyperbolicity of $g$.

If however $ a=b$, this means that the limit set $L (\Gamma)$ 
consists of a single point
$a$ (fixed by each $g \in \Pi$) : contradicting that $\Gamma$ is of the I Kind.

\QED

\section{Mappings to curves }

If a manifold $X$ is fibred ( with connected fibres) onto a curve $C$, certainly 
we have a surjection of fundamental groups $\pi_1 (X) \ra \pi_1 (C)$.

In fact, let $p_1, \dots  p_r$ be the critical values
of $f$, and set $ C^* : = C - \{p_1, \dots  p_r \}$, $ X^* : = f^{-1}( C^* )$:
then we have  surjections  $\pi_1 (X^*) \ra \pi_1 (X)$ ,  $\pi_1 (C^*) \ra \pi_1 (C)$,
 and an exact homotopy sequence 
$$ \pi_1(F) \ra \pi_1(X^*) \ra \pi_1(C^*)\ra 1.$$
It suffices to observe that  the surjection  $\pi_1 (X^*) \ra \pi_1 (C)$ 
factors through $\pi_1 (X) \ra \pi_1 (C)$.

If however there is a non trivial holomorphic map $ F : C \ra C'$, such that $F$ 
does not factor as $ F = F' \circ F"$, where $F'$ is unramified, then it is not
difficult to show that there is a surjection of fundamental groups 
$F_* : \pi_1 (C) \ra \pi_1 (C')$.

This is the reason why one needs some extra assumptions on a surjection of
 fundamental groups $\pi_1 (X) \ra \pi_1 (C')$ in order to decide whether 
the corresponding map is a fibration (i.e., it has connected fibres).

Recall (cf. e.g.  \cite{cko}, Lemma 3 , page 283) the following

\begin{df}
 Let $ m_i$ ,$ \forall i = 1, \dots  r,$ be the G.C.D. of the multiplicities 
of the components of the divisor $ f^{-1} (p_i)$ ($p_1, \dots  p_r$ 
are again the critical values of $f$).

Then the orbifold fundamental group 
$ \pi_1^{orb}(f) $ is defined as the quotient of $\pi_1 ( C - \{p_1, \dots  p_r \})$ 
by the subgroup normally generated by $ \{ \gamma_i^{m_i} \}$, $ \gamma_i$ being a 
simple geometric path around the point $p_i$. 
\end{df}

As a corollary of the results of the previous section we have

\begin{lem}
If  $X$ admits a
surjective holomorphic map $f$ with connected fibres $f: X \rightarrow C$ where C
is a Riemann surface of genus $g \geq 0$, then 
the induced homomorphism $f_* :
\pi_1(X) \ra \Pi_{g}$ is surjective, and   its kernel $H$ 
is finitely generated exactly when $g=0$ or  when 
$g \geq 1$ and there are no multiple fibres, i.e., $ \pi_1^{orb}(f) \cong \pi_1(C)$. 
\end{lem}

\Proof
As well known (cf. e.g. a slightly general version given \cite{cko}, Lemma 3 , page 283,
whose notation we will follow)  
we have an exact sequence 
$$ \pi_1(F) \ra \pi_1(X) \ra \pi_1^{orb}(f) \ra 1$$
where $F$ is a smooth fibre of $f$, and, $p_1, \dots  p_r$ being the critical values
of $f$, and $ m_i, \forall i = 1, \dots  r,$ being the respective G.C.D. of the multiplicities 
of the components of the divisor $ f^{-1} (p_i)$ ,the orbifold fundamental group 
$ \pi_1^{orb}(f) $ is defined as the quotient of $\pi_1 ( C - \{p_1, \dots  p_r \})$ 
by the subgroup normally generated by $ \{ \gamma_i^{m_i} \}$, $ \gamma_i$ being a 
simple geometric path around the point $p_i$. 

Let $\psi : = f_*$.
 Thus $ ker (\psi)$ contains  the normal subgroup $K$, image of $\pi_1(F)$, 
which is  finitely generated since  $F$ is compact, and the cokernel  
$ ker (\psi)/  K $ is isomorphic to the kernel of 
$\rho : \pi_1^{orb}(f) \ra \pi_1(C)$.

Therefore $ ker (\psi)$ is finitely generated if and only if
 $ ker \rho$ is finitely generated. This is then the case for $g = 0$, so
let us assume that $ g\geq 1$.
If we moreover assume that there are no multiple fibres,
then $\rho$ is an isomorphism, and we are again done.

Otherwise, $\pi_1^{orb}(f) $ is a Fuchsian group and the same proof as in Lemma 3.4
shows that $\pi_1^{orb}(f) $  is $ NINF $. Then $ ker \rho$ is finitely generated
if only if it is trivial: since the alternative that its index be finite is 
ruled out by the condition $g\geq 1$.

Finally, if $ ker \rho$  is trivial, then $ \rho$ is an isomorphism,
and by looking at the Abelianization we see that $r =1$.
But then the orbifold fundamental group  $\pi_1^{orb}(f) $ is a free product of
a free group of rank $2g-1$ with a cyclic group of order $m_1$ and its Abelianization
is then not a free Abelian group of rank $2g$, a contradiction.

\QED

\begin{teo}
Let $X$ be a compact K\"ahler manifold, and assume that its fundamental
group admits a non trivial homomorphism $\psi$ 
to the fundamental group $\Pi_g$ of a compact
 Riemann surface of genus $g \geq 2$, with kernel $H$.

Then the following conditions are equivalent 
\begin{itemize}
\item
1) $\psi$ is induced by an irrational pencil of genus $g$  
 i.e. there is a
surjective holomorphic map $f$ with connected fibres 
$f: X \rightarrow C$ such that  $\psi = f_*$, and where C
is a Riemann surface of genus $g$.

\item
2) $\psi$ is surjective and the image of $\psi^* : 
H^1 (\Pi_g , \Q) \ra H^1 (X , \Q)$
contains a $g$-dimensional maximal isotropic subspace
(for the bilinear pairing $H^1 (X , \Q) \times H^1 (X , \Q) \ra 
H^2 (X , \Q)$)
\item
3) $\psi$ induces an injective map in cohomology $\psi^* : 
H^1 (\Pi_g , \Q) \ra H^1 (X , \Q)$, and the image of $\psi^* $ 
contains a $g$-dimensional maximal isotropic subspace

\end{itemize}
likewise,  the following conditions are also equivalent to each other
\begin{itemize}
\item
4) $\psi$ is induced by an irrational pencil of genus $g$  without multiple fibres,
 i.e. ,
 for each fibre $F'$
the equation of divisors $ F' = r D$, with $ r \geq 1$, implies $ r=1$.

\item
5) $\psi$ is surjective and its kernel $H$ is finitely generated

\end{itemize}

\end{teo}

\Proof
We observe first that 2) implies 3) follows since a surjective
homomorphism induces a surjective homomorphism between the abelianizations,
and dualizing one obtains an injective homomorphism in cohomology.

Recall that, by the isotropic subspace theorem of \cite{cat1}, 
given a maximal isotropic 
subspace $V \subset H^1 (X , \Q)$, of dimension $g$, there is a holomorphic fibration
onto a curve $C$ of genus $g$ such that $V \subset f^* (H^1 (C , \Q))$.

Now,  1) implies 2) because, if the pull back of a maximal isotropic
subspace from the given curve $C$ is not maximal, then we have another fibration
$ f' : X \ra C'$ to a curve of genus $g'  > g$ such that $f^* (H^0 (C , \Omega^1_C)) 
\subset f^* (H^0 (C' , \Omega^1_{C'})$, whence $f$ factors through $f'$, contradicting
the fact that $f$ has connected fibres.

Let us show that 3) implies 1). The isotropic subspace theorem gives us 
the desired $ f : X \ra C$, where $C$ has genus $g$.

Since however $C$ is a classifying space for $\Pi_g$, there is a continuous
map $F : X \ra C$ such that $\psi = F_*$.

Compose both maps with the Jacobian embedding  $ \alpha : C \ra J $ and observe 
that by the proof 
of the isotropic subspace theorem the two subspaces $ f^* (H^1 (C , \Q))$ and
 $ F^* (H^1 (C , \Q))$ coincide. 

Therefore the two maps $ \alpha \circ F$,$ \alpha \circ f$ are given by integrals 
of the same differentiable 1-forms,
hence there is an isogeny $ p : J \ra J$ such that $ \alpha \circ f = p \circ 
\alpha \circ F$. We get thus, up to changing $F$ in its homotopy equivalence,
a factorization  $f = p' \circ F$.  Thus, $p'$ is surjective, and actually 
it has degree $1$ or otherwise $ f^* (\eta)$ , with $(\eta)$ the positive 
generator of $H^2 (C , \Z))$, would be divisible, contradicting that 
$f$ has connected fibres.

The conclusion is that $f$ and $F$ are homotopy equivalent, thus 
$\psi = F_*= f_*$.

The implication 4) $\ra$ 5) is exactly lemma 4.2, so we are left with
 showing that 5) implies 4).

Now, 5) implies that the image of $\psi^*$ contains a $g$-dimensional 
isotropic subspace. Assume this subspace not to be maximal: then there is a
fibration $f : X \ra C$ where the genus $g'$ of $C$ is strivtly larger than $g$.

Argueing as we did before, we find a factorization of $\psi$ through $f'_*$.

Since $\psi$ is surjective, 
Lemma 3.2. applies and we get that $g' = g$, and $\psi = f_*$.

Finally, $f$ has no multiple fibres again by Lemma 4.2.

\QED

\section{The logarithmic case}

In this section we shall generalize the results of the previous
section to the case where we have a Zariski open set $Y$ 
in a compact K\"ahler manifold $X$. One may assume w.l.o.g.
that the complement $ X - Y$ is a normal crossings divisor
$ D$. We shall consider holomorphic maps $ f : Y \ra C$,
where $C$ is  Zariski open in a compact curve $\bar{C}$, and
the map $f$ is meromorphic on $X$, whence there is another
 compactification $\bar{X}$ of $X$ where $f$ extends
 holomorphically.

When we shall say that $f$ is a {\bf pencil}, we shall mean that 
$f$ is as above, that the extension $\bar{f}$ of $f$ has connected fibres, and that
$f$ is surjective (Arapura calls these maps { \bf admissible maps}).
We shall denote by $B$ the complement $\bar{C} - C$,
because quite often it will be the branch locus of
a fibration of a compact manifold. 

However, $X$ will not necessarily be non compact, 
the reason for this being that we shall here consider
surjective homomorphisms $ \pi_1(X) \ra \F_n$ to a non 
Abelian free group.

Notice moreover that 
\begin{itemize}
\item
any automorphism of $\Pi_g$ composed with
the standard surjection $p : \Pi_g \ra \F_g$ such that $ p(a_i) =
p (b_i) = x_i$ produces a maximal isotropic subspace
of $ H^1(C, \Z)$.
\item
there is a surjection  $p : \Pi_g \ra \F_n$ iff $g \geq n$ (since
$Im (p^*)$ is an isotropic subspace of dimension $n$).

\end{itemize}

The next theorem extends the results of I. Bauer and D. Arapura 
(Thms 2.1. and 3.1 of \cite{ingrid} and Cor. 1.8 of \cite{ara}, cf. also
thm. 2.11 of \cite{cat3}) using the new ideas introduced in the previous 
sections. 

Observe that also in this context a pencil induces a surjective 
homomorphism of fundamental groups.

\begin{df}
Let $f : Y \ra C$ be a pencil as above. We shall say that $f$ is 
of type $g$ if either 
\begin{itemize}
\item
$C$ is compact of genus $g$ ( then $\pi_1(C) \cong \Pi_g$ )
\item
$C$ is not compact and its first Betti number equals $g$ 
( then $\pi_1(C) \cong \F_g$ ).

\end{itemize}

\end{df}

\begin{df}
Let $f : Y \ra C$ be a pencil as above. 

We may assume that $f$ extends to a holomorphic fibration
$F : X \ra \bar{C}$. We can separate the complementary divisor
$D = X - Y$ into three parts :

\begin{itemize}
\item
$D^{hor}$ : the union of the components dominating $C$
\item
$D^{FVert}$: the union of the fibres over $\bar{C} - C$
\item
$D^{PVert}$: the union of the components mapping to points of $C$.
\end{itemize}

We define then the orbifold fundamental group  of $f$  by the usual
procedure:
let $p_1, \dots  p_r$ be the critical values
of $f$, and $ m_i, \forall i = 1, \dots  r,$ be the respective G.C.D. 
of the multiplicities 
of the components of the divisor $ f^{-1} (p_i)$.
Then the orbifold fundamental group 
$ \pi_1^{orb}(f) $ is defined as the quotient of 
$\pi_1 ( C - \{p_1, \dots  p_r \})$ 
by the subgroup normally generated by $ \{ \gamma_i^{m_i} \}$, 
$ \gamma_i$ being a 
simple geometric path around the point $p_i$. 

\end{df}

\begin{lem}
Let $X$ be a compact K\"ahler manifold and $Y = X - D$ be a Zariski open 
set. If  $Y$ admits a
pencil  $f: Y \rightarrow C$, then 
the induced homomorphism $f_* :
\pi_1(X) \ra \pi_1(C)$ is surjective, and   its kernel $H$ 
is finitely generated exactly when $g=0$ or  when 
$g \geq 1$ and there are no multiple fibres, i.e., 
$ \pi_1^{orb}(f) \cong \pi_1(C)$. 
\end{lem}

\Proof
We use once more the exact sequence (the situation being more general,
but the proof exactly the same as in \cite{cko}, Lemma 3 , page 283)
$$ \pi_1(F) \ra \pi_1(X) \ra \pi_1^{orb}(f) \ra 1$$
where $F$ is a smooth fibre of $f$ which is transversal to $D^{hor}$.

Let $\psi : = f_*$.
 Thus $ ker (\psi)$ contains  the normal subgroup $K$, image of $\pi_1(F)$, 
which is  finitely generated since $F$ is of finite type, and the cokernel  
$ ker (\psi)/ K $ is isomorphic to the kernel of 
$\rho : \pi_1^{orb}(f) \ra \pi_1(C)$.

Therefore $ ker (\psi)$ is finitely generated if and only if
 $ ker \rho$ is finitely generated. 
Assume that there are no multiple fibres:
then $\rho$ is an isomorphism.

Otherwise, $\pi_1^{orb}(f) $ is a Fuchsian group and the 
same proof as in Lemma 3.4
shows that $\pi_1^{orb}(f) $  is $ NINF $. Then $ ker \rho$ is finitely generated
if only if it is trivial,  its index being infinite for $g \geq 1$.

\QED

\begin{teo}
Let $X$ be a compact K\"ahler manifold and $Y = X - D$ be a Zariski open 
set. Assume that the fundamental group of $Y$  admits a 
homomorphism $\psi : \pi_1(Y) \ra \F_g$ to a free group of rank $g$, 
with kernel $H$.

Then the following are equivalent 

\begin{itemize}
\item
1) $\psi$ is induced by a pencil $ f : Y \ra C$ of 
type $g$ and by a surjection $\pi_1(C) \ra \F_g$.
\item
2) $\psi$ is surjective and the image of $\psi^* : 
H^1 (\F_g , \Q) \ra H^1 (Y , \Q)$
is a $g$-dimensional maximal isotropic subspace
(for the bilinear pairing $H^1 (Y , \Q) \times H^1 (Y , \Q) \ra 
H^2 (Y , \Q)$)

\end{itemize}
likewise,  the following are also equivalent to each other
\begin{itemize}
\item
1') $\psi$ is induced by a pencil of type $g$  
without multiple fibres

\item
2') $\psi$ is surjective and the kernel $H$ of $\psi$ 
is finitely generated

\end{itemize}
Finally, the curve $C$ is compact if and only if 
$\psi^* ( 
H^1 (\F_g , \Q) )$ is  also an isotropic subspace  in $ H^1 (X , \Q) $ $(
\subset  H^1 (Y , \Q) )$.

\end{teo}

\Proof

1) implies 2) : it suffices to show that $\psi^* ( H^1 (\F_g , \Q) )$ 
is a maximal isotropic subspace. Assume the contrary : then, there 
is a strictly larger maximal isotropic subspace $V$ induced
(cf.\cite{cat3}, thm. 2.11) by a pencil $f': Y \ra C'$.

The pencil is induced by integrations of linearly independent forms 
in $H^0(\Omega^1( log D)$, whence we get a factorization
$ Y \ra C' \ra C$: since $f$ has connected fibres $ C' \cong C$,
contradicting that the logarithmic genus of $C'$ is strictly larger 
than $g$.

Assume 2): then by Arapura's cor. 1.9 there is a pencil $f: Y \ra C$
and $ \tau : H_1(C, \Z) \ra \Z^g$ such that there is a 
factorization in homology 
$H_1(\psi) = \tau \circ H_1(f)$.

Since $\psi$ is surjective, then also $H_1(\psi)$ is such, whence we
may lift $\tau$ to a surjection $ \pi_1(C) \ra \F^g$ (since $Aut(\F_g)$
surjects onto $GL (g, \Z)$, cf. \cite{mks}, sections 3.5 and 3.6).

The equivalence of 1') and 2') follows as in the proof of
Theorem 4.3 in view of Lemma 5.3 .
The last assertion is already contained in theorem. 2.11, loc. cit.

\QED

\bigskip

\section{Fibred algebraic surfaces and  good open sets}

In this section we shall consider a smooth compact algebraic surface,
and a holomorphic fibration $ f : S \ra C$.

\begin{df}
A {\bf good open set} of a fibration will be any set of the form:
$ U = f^{-1} (C-B)$, where $B$ is any finite set containing the
set of critical values of $f$.
\end{df}

\begin{oss}
In the above situation one has an exact sequence of fundamental groups
$$ 1 \ra \pi_1(F) \ra \pi_1(U) \ra \pi_1(C-B) \ra 1$$
where $F$ is any fibre of $f$ over a point of $C-B$.
\end{oss}

The next theorem will give a topological characterization of 
good open sets of some fibration.

The case where $U= S$ was already treated by Kotschick 
(\cite{kot}, Proposition 1) and Hillman (\cite{hil}), cf. also Kapovich (\cite{kap}),
and we only prove   an $\epsilon^2$ more general result:

\begin{teo}
Assume that $S$ is a compact K\"ahler surface and that 
we have an exact sequence  
$$ 1 \ra \Pi_r \ra \pi_1(S) \ra \Pi_g \ra 1$$
where $g \geq 2$. If moreover the topological 
Euler-Poincare' characteristic of $S$,
$e(S)$, equals $ 4 (g-1) (r-1)$,
then there exists a holomorphic submersion $f: S \ra C$ inducing 
the previous exact sequence.
\end{teo}

\Proof
By Theorem 4.3 we find a fibration $ f : S \ra C$, where $C$ has genus $g$,
inducing the given epimorphism of fundamental groups.

By the theorem of Zeuthen-Segre the Euler-Poincare' characteristic of $S$,
$e(S)$, equals $ 4 (g-1) (s-1) + \mu$ where $s$ is
the genus of a smooth fibre of $f$ and where $\mu \geq 0$, equality holding
if and only if all the singular fibres are multiple of 
a smooth elliptic curve. It is clear that $ s \geq r$.

Our assumption implies $ 4 (g-1) (s-r) + \mu = 0$, whence $s=r$ and
$\mu =0$. We are therefore done in the case where $r=0$ or $r\geq 2$.

If finally $r=1$, we are done unless there is some multiple fibre.
But the existence of a multiple fibre is excluded by Lemma 4.2.

\QED

We come now to the non compact case

\begin{teo}
Assume that $U$ is a non complete Zariski open set 
of an algebraic surface
 and that the following properties hold:

\begin{itemize}
\item
(P1) we have an exact sequence  
$$ 1 \ra \Pi_r \ra \pi_1(U) \ra \F_g \ra 1$$
where $g \geq 2$. 
\item
(P2) The topological 
Euler-Poincare' characteristic of $U$,
$e(U)$, equals $ 2 (g-1) (r-1)$.
\item
(P3) For each end $\EE$ of $U$, the corresponding  fundamental group 
$\pi_1^{\EE}$ surjects onto a cyclic subgroup of $\F_g$, and each
simple geometric generator $\gamma_i$ has a non trivial image
in $\F_g$.

\end{itemize}
Then $U$ is a good open set of a fibration, more precisely,
there exists a proper holomorphic submersion $f: U \ra C$ inducing 
the previous exact sequence.
\end{teo}

\Proof
By Thm. 5.4 we have a fibration $f : U \ra C$ inducing the surjection
$\pi_1(U) \ra \F_g \ra 1$, and w.l.o.g. we have an extension
$\bar{f} : \bar{S} \ra \bar{C}$, where $ \bar{S}$ is a blow-up of $S$.  
By condition (P3) there is no component of $D$ which is horizontal,
and each component of $D$ pulls back to a fibre of $\bar{f}$.

It turns therefore out that there is no point of indeterminacy of $f$
on $D$, whence $ \bar{S} = S$, and again by condition (P3) $U$ is the 
full inverse image of $C$ under $\bar{f}$.

It suffices to apply the logarithmic version of the Zeuthen Segre
theorem, similarly to thm. 2.14 of \cite{cat3}, and we conclude that
( $e(U)$ being the same in ordinary and Borel-Moore homology by 
virtue of Poincare'
duality ) $e(U) = 2 (g-1) (r-1) \geq 2 (g-1) (s-1) +\mu$, where
$s \geq r$ is the genus of a smooth fibre of $f$. Whence, as usual,
$0 \geq 2 (g-1) (s-r) +\mu$, thus $s=r$ and $\mu = 0$.
We conclude as in theorem 6.3.

\QED

\bigskip

\begin{footnotesize}
\noindent

{\bf Note.}
The next question : when is the fibration $f$ 
a constant moduli fibration?
was already answered, with similar methods, in the previous paper \cite{cat3}, cf.
5.4 and 5.7.

\end{footnotesize}
 
\section{Restrictions for the monodromy}

Before we present some interesting corollary of the previous theorem,
 we need to recall some well known results

\begin{lem}
Let $X$ be a topological manifold and $\Gamma$ its fundamental group.
Then $H^1(X, \Z) = H^1(\Gamma, \Z)$ and $H^2(\Gamma, \Z)$ injects into
$H^2(X, \Z) .$
\end{lem}

\Proof
Let $\tilde{X}$ be the universal covering of $X$, so that 
$ X \cong \tilde{X}/ \Gamma$.

The proof is a direct consequence of the spectral sequence for group cohomology
with terms $ H^p (\Gamma, H^p (\tilde{X}, \Z))$, converging to a suitable graded quotient of 
$ H^{p+q}(X, \Z)$, in view of the fact that $H^1 (\tilde{X}, \Z)= 0$.

\qed

\begin{lem}
Let $ 1 \ra H \ra \Gamma \ra B \ra 1$ be an exact sequence of groups, where $B$ 
is a finitely generated free group, 
or the fundamental group $\Pi_{g}$ of
 a compact hyperbolic Riemann surface.
Assume that:

 (**) $ H^2(B, \Z)$ injects into $H^2(\Gamma, \Z)$.

Then $H^1(B, \Z) \subset H^1(\Gamma, \Z)$ with quotient $ H^1(H, \Z)^B$, and 
the cup product $ H^1(B, \Z)\times H^1(\Gamma, \Z) \ra  H^2(\Gamma, \Z)$ 
lands into 
 the subgroup $F$ fitting into the exact sequence
 $ 0 \ra H^2(B, \Z) \ra F \ra H^1(B, H^1(H, \Z)) \ra 0.$

In particular, 
 let $V$ a maximal isotropic subspace of $ H^1(B, \Z)$: then
$V$ remains a maximal isotropic subspace in $ H^1(\Gamma, \Z)$

only if (resp. : if and only if, in the case where $B$ is free)

(***) the cup product 
$ H^1(B, \Z) \times H^1(H, \Z)^B \ra H^1(B, H^1(H, \Z))$ is non degenerate
in the second factor.


\end{lem}

\Proof

The proof of the first assertions is a direct consequence of the spectral 
sequence for group cohomology, with terms
$ H^p (B, H^q (H, \Z))$, converging to a suitable graded quotient of 
$ H^{p+q}(\Gamma, \Z)$, in view of the fact that
 , by assumption (**), the differential 
 $d_2 : H^1(H, \Z)^B= H^0 (B, H^1 (H, \Z))  \ra H^2 (B, \Z)$ is zero.

The second assertion holds with "if and only if" in the case where
 $B$ is a free group, since then 
 $ H^2(B, \Z) = 0$, $ F = H^1(B, H^1(H, \Z))$, and the question is whether
 $ H^1(B, \Z)$ is a maximal isotropic subspace in $ H^1(\Gamma, \Z)$.

In the other case, observe that  $ H^2(B, \Z) = \Z$, and that we can 
find two maximal isotropic subspaces $V, V'$ such that
  $ H^1(B, \Z) =  V \bigoplus  V' $ : moreover then the cup 
product yields an isomorphism of $V'$ with  $V^{\vee}$. 

If we get an element $w' \in H^1(H, \Z)^B $ annihilating $ H^1(B, \Z)$, this means
that there is a lift $ w \in H^1(\Gamma, \Z) - H^1(B, \Z)$ such that $ w \cup H^1(B, \Z)
\subset H^2(B, \Z)$.
In particular, there is $u \in V^{\vee}$ such that $(w-u) \cup V = 0$.

We easily conclude then that the span of $V$ and of $w-u$ is isotropic.


\qed

\begin{cor}
If the finitely presented group $\Gamma$ admits a surjection 
$\Gamma \rightarrow 
\Pi_{g}$ with finitely generated kernel $H$, then $\Gamma$ cannot be 
the fundamental group of a compact K\"ahler manifold $X$ if there is
a non zero element $ u \in  H^1(H, \Z)^{\Pi _g }$ such that cup product 
with $u$ yields the identically 
  zero map $$  H^1(\Pi _g , \Z) \ra H^1(\Pi _g , H^1(H, \Z)).$$

\end{cor}

We want now to explicitly write down, for $\Pi$ equal
either to $\Pi_g$ or to a free group $\F_g$,  the condition that there is
a non zero element $ u \in  H^1(H, \Z)^{\Pi }$ such that cup product 
with $u$ yields the identically 
  zero map 
$$  H^1(\Pi , \Z) \ra H^1(\Pi , H^1(H, \Z)).$$

Observe first that $H^1(\Pi , \Z) = Hom_{\Z} (\Pi, \Z) \cong \Z^b$,
where $b=g$ in the free case, otherwise $b=2g$.

The condition that, $\forall \phi \in Hom_{\Z} (\Pi, \Z)$, $ \phi u =0$
 in $H^1(\Pi , H^1(H, \Z))$
means that there is an element $v_{\phi} \in H^1(H, \Z))$ such that
$$ \forall \gamma \in \Pi, \  \phi (\gamma) u = \gamma v_{\phi} - v_{\phi}.$$

Taking a basis $\phi_1, \dots \phi_b$, we get $v_1, \dots v_b$ such that
$$ (1) \  \forall \gamma \in \Pi, \   \gamma v_j = v_j + \phi_j (\gamma) u.$$
Recall moreover that $u$ is invariant, whence 
$$ (2) \ \forall \gamma \in \Pi, \   \gamma u =  u.$$
Conditions (1), (2), and the $\Z$-linear independence of the characters
$\phi_j$ implies the $\Z$-linear independence of $u, v_1, \dots v_b$ since
the $\Z$-module $H^1(H, \Z))$ is torsion free.

\begin{df}
A bad monodromy module is a free $\Z$-module of rank $b+1$, with basis
$u, v_1, \dots v_b$, and with an action of $\Pi$ given by (1) and (2).
\end{df}

\begin{ex}
Let $H$ be a finitely generated group, and let 
$ 1 \ra H \ra \Gamma \ra \Pi_g \ra 1$ be  an exact sequence such that 
the induced  action of $\Pi_g$ on $H$ by conjugation induces on the $\Z$ 
dual of the Abelianization of $H$ a $\Pi_g$ -module structure which
contains a bad monodromy module.
Then $\Gamma$ cannot be the fundamental group of a compact K\"ahler manifold.
\end{ex}

\begin{oss}
One can use the same type of restriction in the case where $U \neq X$ is the
inverse image of the non critical values of a fibration $f : X \ra \bar{C}$,
 and obtain 
in this way a restriction for the monodromy in the case where $\bar{C}-B$
has first Betti number at least $2$. 
\end{oss}

\begin{oss}
To see finally the relation of the above condition with the
 theory of Lefschetz pencils (in particular with the splitting 
in invariant and vanishing cycles),
let us observe that our cup product is non degenerate if we have a monodromy
invariant splitting $ H^1(H, \Z) = H^1(H, \Z)^{\Pi } \bigoplus W$.

Because then we may write $ v_{\phi} = u_{\phi} + w_{\phi} $ and we obtain
$$ \forall \gamma \in \Pi, \  \phi (\gamma) u = \gamma v_{\phi} - v_{\phi} = 
\gamma w_{\phi} - w_{\phi}  \in W  ,$$

whence $ u=0$.

This splitting is proven by Deligne's Semi-simplicity Theorem (4.2.6 of
\cite{del71}, cf. also thm. 9.1, page 37 , Chapter II of \cite{grif84}).
 \end{oss}

\vfill

\noindent
{\bf Author's address:}

\bigskip

\noindent 
Prof. Fabrizio Catanese\\
Lehrstuhl Mathematik VIII\\
Universit\"at Bayreuth, NWII\\
 D-95440 Bayreuth, Germany

e-mail: Fabrizio.Catanese@uni-bayreuth.de

\end{document}